\documentclass[a4paper,10pt]{article}
\usepackage{amsmath}
\usepackage{amssymb}
\usepackage{amsthm}
\usepackage{graphicx}
\usepackage{natbib}
\usepackage{cancel}
\usepackage[affil-it]{authblk}
\def\ci{\perp\!\!\!\perp}

\newtheorem{mydef}{Definition}
\newtheorem{mythm}{Theorem}
\newtheorem{mylemma}{Lemma}

\title{
On the graph-theoretical interpretation of Pearson correlations in a multivariate process and a novel partial correlation measure
}

\author{Jakob Runge}
\affil[1]{Potsdam Institute for Climate Impact Research,
14473 Potsdam, and
Department of Physics, Humboldt University, 12489 Berlin, Germany}

\begin{document}
\maketitle

\begin{abstract}
The dependencies of the lagged (Pearson) correlation function on the coefficients of multivariate autoregressive models are interpreted in the framework of time series graphs. Time series graphs are related to the concept of Granger causality and encode the conditional independence structure of a multivariate process. The authors show that the complex dependencies of the Pearson correlation coefficient complicate an interpretation and propose a novel partial correlation measure with a straightforward graph-theoretical interpretation.
The novel measure has the additional advantage that its sampling distribution is not affected by serial dependencies like that of the Pearson correlation coefficient. In an application to climatological time series the potential of the novel measure is demonstrated.
\end{abstract}

\section{Introduction}
Among the measures of association, the Pearson (product-moment) correlation coefficient is widely applied in many fields of science due its simple computation and alleged ease of interpretation. 
Indeed, the square of this correlation coefficient between two processes simply represents the proportion of variance of one process that can be linearly represented by the other \citep{Chatfield2003}. But what does this value say about how strong both processes are associated or dependent with each other in a multivariate process?
While it is a commonplace that correlation does not imply causation \citep{Spirtes2000}, the aim of this article is to further elucidate how the value of the lagged Pearson correlation coefficient -- in the following referred to as the correlation (function) -- between two causally dependent components of a multivariate process is to be interpreted.

Graphical models \citep{lauritzen1996graphical} provide a well interpretable framework to study interactions in a multivariate process. Here we utilise the derived concept of time series graphs \citep{Dahlhaus2000,Eichler2011} to study the dependencies of cross correlation for the class of multivariate autoregressive time series models in a graph-theoretical way. We demonstrate that cross correlation can be rather misguiding as a measure of how \emph{strong} two processes interact and is ambiguously influenced by other dependencies in the multivariate process.

Based on the time series graph a certain partial correlation measure is introduced for which we prove very simple dependencies on the autoregressive coefficients, making it straightforward to be interpreted as the strength of dependence between these two components alone. We also introduce further partial correlation measures that capture different aspects of the dependence between two components.

Another commonly known problem of cross correlation is the estimation of its significance in the presence of strong autocorrelations in the time series. These dependencies violate the assumption of independent identically distributed samples and `inflate' the sampling distribution making an assessment of significance difficult. For the proposed partial correlation measure, on the other hand, we show analytically and numerically that the is not affected by autocorrelation, as our theoretical results suggest.

The article is structured as follows: In Sect.~\ref{sect:anamod} we define time series graphs and their relation to autoregressive models. In Sect.~\ref{cc_deps} the dependencies of the lagged correlation function are interpreted graph-theoretically. In Sect.~\ref{par_corr} the novel partial correlation measures is introduced and some theoretical results are discussed. The properties of its sampling distribution are investigated in Sect.~\ref{numerics}. Finally, in Sect.~\ref{application} we compare the differences between the measures on a climatological example of temperature time series in the tropics.


\section{Time series graphs and autoregressive models} \label{sect:anamod}

\subsection{Time series graphs}
\begin{figure}[!t]
\begin{center}
\includegraphics[width=.8\columnwidth]{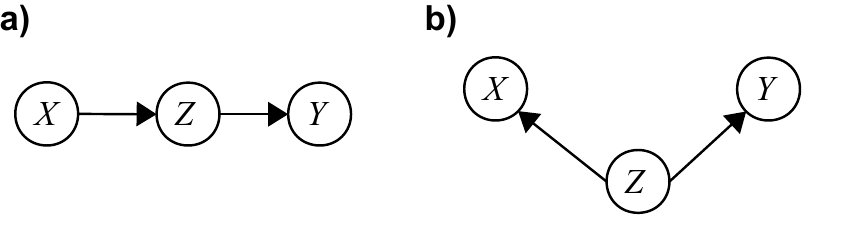}
\end{center}
\caption[]{Causality between three processes: (a) Indirect chain and (b) common driver system.}
\label{fig:causality}
\end{figure}

Graphical models \citep{lauritzen1996graphical} provide a tool to distinguish direct from indirect interactions between and within multiple processes.
Underlying is the concept of \emph{conditional independencies} in a general multivariate process, which can be explained as follows. Consider three processes where $X$ drives $Z$ (i.e., $Z$ is statistically dependent on $X$ at some lag in the past) and $Z$ drives $Y$ as visualised in Fig.~\ref{fig:causality}(a). Here $X$ and $Y$ are not directly but indirectly interacting and in a bivariate analysis $X$ and $Y$ would be found to be dependent -- implying that their correlation coefficient would be non-zero in the case of a linear dependency. The same holds for a common driver scheme in Fig.~\ref{fig:causality}(b). 
If, however, the variable $Z$ is included into the analysis, one finds that $X$ and $Y$ are independent \emph{conditional} on $Z$, written as
\begin{align*}
 X~~\ci~~ Y~~ |~~ Z\,.
\end{align*}

This concept is now applied to define links in a time series graph \citep{Eichler2011} of a multivariate stationary discrete-time process $\mathbf{X}$. Each node in that graph represents a single random variable, i.e., a subprocess, at a certain time $t$. Nodes $X_{t-\tau}$ and $Y_t$ are connected by a directed link ``$X_{t-\tau} \to Y_t$'' pointing forward in time if and only if $\tau>0$ and 
\begin{align}  \label{eq:def_graph}
  X_{t-\tau}~~ \cancel{\ci}~~ Y_t~~ |~~ \mathbf{X}_t^-\setminus \{X_{t-\tau}\},
\end{align}
i.e., if they are \emph{not} independent conditionally on the past of the whole process denoted by $\mathbf{X}^-_t=(\mathbf{X}_{t-1},\,\mathbf{X}_{t-2},\ldots)$.
If $Y\neq X$, the link ``$X_{t-\tau} \to Y_t$'' represents a \textit{coupling at lag} $\tau$, while for $Y=X$ it represents an \textit{autodependency at lag} $\tau$. 
Further, nodes $X_t$ and $Y_t$ are connected by an undirected contemporaneous link ``$X_{t} \-- Y_t$'' \citep{Eichler2011} if and only if
\begin{align} \label{eq:def_graph_contemp}
X_t~~ \cancel{\ci}~~ Y_t~~ |~~ \mathbf{X}_{t+1}^- {\setminus} \{X_t,Y_t\},
\end{align}
where also the contemporaneous present $\mathbf{X}_t{\setminus}\{X_t,Y_t\}$ is included in the condition. 
Note that for stationary processes it holds that ``$X_{t-\tau} \to Y_t$'' whenever ``$X_{t'-\tau}\to Y_{t'}$'' for any $t'$.

These graphs can be linked to the concept of a lag-specific \emph{Granger causality} \citep{Granger1969,Eichler2011,Runge2012prl}.
In the original definition of Granger causality $X\in \mathbf{X}$ \emph{Granger causes} $Y\in \mathbf{X}$ with respect to the past of the whole process $\mathbf{X}$ if (1) events in $X$ occur before events in $Y$ and (2) $X$ improves forecasting $Y$ even if the past of the remaining process $\mathbf{X}{\setminus}\{X,Y\}$ is known. The latter property is directly related to the conditional dependence between $X$ at some lag and $Y$ given the past of the remaining process $\mathbf{X}{\setminus}\{X,Y\}$ which defines links in the time series graph.
In \cite{Eichler2005,Eichler2011} the range and conditions of application are further discussed.

For the following analysis the notion of \emph{parents} $\mathcal{P}_{Y_t}$ and \emph{neighbors} $\mathcal{N}_{Y_t}$ of a process $Y_t$ in the time series graph will be important.
They are defined as
\begin{align}
\mathcal{P}_{Y_t} &\equiv \{Z_{t-\tau}:~ Z\in \mathbf{X},~\tau>0,~Z_{t-\tau}\to Y_t\},\\
\mathcal{N}_{Y_t} &\equiv \{X_t:X\in\mathbf{X}, X_t {\--} Y_t\}.
\end{align}
Note, that also the past lags of $Y$ can be part of the parents. 
The parents of all subprocesses in $\mathbf{X}$ together with the contemporaneous links comprise the time series graph.

\subsection{Relation to multivariate autoregressive models}

While the definition of time series graphs was given for the large class of processes sufficing condition (S) in \cite{Eichler2011}, in this article we consider the case of a stationary $N$-variate discrete-time process defined as
\begin{align} \label{eq:var}
\mathbf{X}_t &= \sum_{s=1}^p \Phi(s) \mathbf{X}_{t-s} + \varepsilon_t ~~~~~~~\varepsilon_t\sim \mathcal{N}(0,\Sigma), 
\end{align}
i.e., a vector autoregressive process of order $p$ where $\Phi(s)$ are $N \times N$ matrices of coefficients for each lag $s$ and the $N$-vector $\varepsilon$ is an independent identically distributed Gaussian random variable with zero mean and covariance matrix $\Sigma$. $\mathbf{\varepsilon}$ is sometimes referred to as the \emph{innovation term}. Its variances on the main diagonal of $\Sigma$ we denote by $\sigma_{i}^2$ and the covariances by $\sigma_{ij}$ for $i\neq j$. 

For this model class the directed and contemporaneous links of the corresponding time series graph are defined by non-zero entries in the coefficient matrix $\Phi$ and the inverse of the innovation covariance matrix $\Sigma$ \citep{Eichler2011}:
\begin{align} \label{eq:def_graph_ar}
X_{t-\tau} \to Y_t ~~~~&\Leftrightarrow~~~~ \Phi_{YX}(\tau) \neq 0 \\
X_{t} \-- Y_t ~~~~&\Leftrightarrow~~~~ (\Sigma^{-1})_{YX} \neq 0.
\end{align}
An alternative definition of contemporaneous links is based on non-zero entries in $\Sigma_{YX}$ \citep{Eichler2011}.

As an example, consider the bivariate autoregressive model of order 1
\begin{align} \label{eq:ar_model_matrix}
\left( \begin{matrix} X_t \\ Y_t \end{matrix} \right) &= \underbrace{\left( \begin{matrix} a & 0\\ c & b \end{matrix} \right)}_{\Phi(1)} \left( \begin{matrix} X_{t-1} \\ Y_{t-1} \end{matrix} \right) + \left( \begin{matrix} \varepsilon_{X,t} \\ \varepsilon_{Y,t} \end{matrix} \right)
\end{align}
and $\Phi(s)=0$ for $s>1$.
\begin{figure}[!t]
\begin{center}
\includegraphics[width=.8\columnwidth]{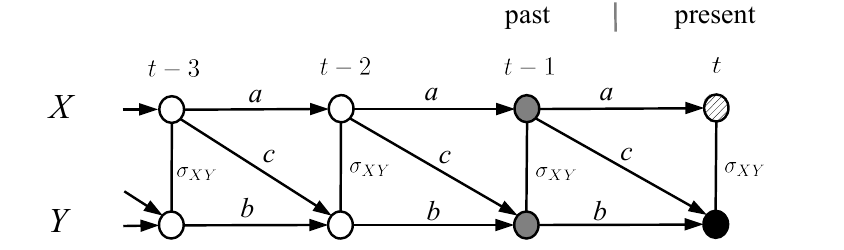}
\end{center}
\caption[]{Visualization of model Eq.~(\ref{eq:ar_model_matrix}) as a time series graph. The labels indicate the coefficients in the matrices $\Phi(1)$ and $\Sigma$. Note, that a non-zero coefficient only determines the existence of absence of a link, but not a weight. Note, that a non-zero $\sigma_{XY}$ only defines a contemporaneous link in the bivariate case, while it is non-zero entries in $(\Sigma^{-1})_{YX}$ in the multivariate case. Due to stationarity, links for $t$ imply links for all $t-1,\,t-2,\,\ldots$. Process $Y_t$ (black node) has one neighbor $X_t$ (hatched node) and two parents (gray nodes).}
\label{fig:tsg_ar}
\end{figure}
In Fig.~\ref{fig:tsg_ar} the corresponding time series graph is visualised. Note, that a non-zero coefficient in the matrices $\Phi_{YX}(\tau)$ or $(\Sigma^{-1})_{YX}$ only defines the existence or absence of a link. In the next sections we address the question of how the \emph{weight} of a link can be quantified. 

\section{Cross correlation of a multivariate autoregressive process} \label{cc_deps}
We are interested in the cross correlation lag function of stationary zero-mean random variables $X,\,Y$ given by
\begin{align} \label{cross_corr}
\rho_{YX}(\tau) \equiv \frac{E[X_{t+\tau}Y_t]}{\sqrt{E[Y_t Y_t]}\sqrt{E[X_t X_t]}},
\end{align}
which depends on the covariances and variances. Thus, we will now give an interpretation of the lagged covariance structure of a multivariate autoregressive process in the framework of time series graphs.

\subsection{Interpretation in terms of paths} \label{sect:cc_deps_paths}
For an autoregressive process given by Eq.~(\ref{eq:var}) there exists an analytical expression of the lagged covariance in terms of $\Phi$ \cite[Ch.~11.3]{brockwell2009time}:
\begin{align} \label{eq:anacov}
\Gamma_{ij}(\tau) \equiv E[\mathbf{X}^i_{t+\tau}  \mathbf{X}^j_t] = \sum^\infty_{n=0} \left(\Psi(n+\tau) \Sigma \Psi^\top(n) \right)_{ij}
\end{align}
where $\Psi(n)$ can be recursively computed from matrix products:
\begin{align} \label{eq:psi_matrices}
\Psi(n) &\equiv \sum_{s=1}^n \Phi(s) \Psi(n-s),
\end{align}
for example,
\begin{align}
\Psi(0) &= \mathbb{I}, \nonumber\\
\Psi(1) &= \Phi(1),\nonumber\\
\Psi(2) &= \Phi^2(1) + \Phi(2),\nonumber\\
\Psi(3) &= \Phi^3(1)  + \Phi(1)\Phi(2) + \Phi(2)\Phi(1)+\Phi(3),
\end{align}
where $\mathbb{I}$ is the identity matrix.

Now, like a non-zero entry in $\Phi$ corresponds to a link, an entry $\Psi(3)_{ij}\neq0$ can be interpreted as a superposition of the contributions from different paths in the time series graph, each with total delay 3: one direct path of only one link with lag 3 [$\Phi(3)_{ij}$], paths composed of two links where the first has lag 1 and the second lag 2 [($\Phi(1)\Phi(2))_{ij}$] and vice versa [($\Phi(2)\Phi(1))_{ij}$], and paths comprised of three links, each with lag 1 [$(\Phi^3(1))_{ij}$]. For example, in the model Eq.~(\ref{eq:ar_model_matrix}), $\Psi(n)$ is given by  $\Phi(1)^n$ and a non-zero coefficient $[\Phi(1)^3]_{YX}\neq0$ thus corresponds to all paths comprised of three links, each with lag 1, e.g., ``$X_{t-3}\to X_{t-2}\to Y_{t-1}\to Y_{t}$''. 
These paths can be interpreted as an indirect causal chain as pictured in Fig.~\ref{fig:causality}(a).

The covariance $\Gamma_{ij}(\tau)$, thus, is an infinite sum of products of $\Psi(n+\tau)$, $\Psi(n)$ and $\Sigma$ and therefore a nonlinear polynomial combination of coefficients of \emph{all possible paths} that end in $\mathbf{X}^j$ and $\tau$-lags later in $\mathbf{X}^i$, emanating from nodes and their contemporaneous neighbors at all past lags. Note, that possible paths via an intermediate node $X_{t-\tau'}$ can only contain the motifs ``$\to~X_{t-\tau'}~\to$'', ``$\--~X_{t-\tau'}~\--$'' or ``$\--~X_{t-\tau'}~\to$'', but not ``$\to~X_{t-\tau'}~\--$'' or ``$\to~X_{t-\tau'}~\leftarrow$'' \citep{Eichler2011}. 

In essence, most non-zero values in the covariance lag function are due to the common driver effect of the past (Fig.~\ref{fig:causality}(b)) or the indirect causality effect due to intermediate lags (Fig.~\ref{fig:causality}(a)). Therefore, the cross correlation as the covariance normalised by the variances, cannot be related to the interaction between $\mathbf{X}^j$ and $\mathbf{X}^i$ alone, i.e., the link ``$\mathbf{X}^j_{t-\tau}\to\mathbf{X}^i_t$'' in the time series graph. Large cross correlation values between two nodes can simply be due to the superposition of indirect paths while the coefficient of the connecting link could be very small (or even zero). In the application (Sect.~\ref{application}) we give an example where this is the case.

\subsection{Interpretation in terms of parents} \label{sect:cc_deps_par}
One can also characterise the dependencies of the covariance Eq.~(\ref{eq:anacov}) in terms of the parents in the time series graph. 

Two univariate subprocesses $X,Y$ of $\mathbf{X}$ given by Eq.~(\ref{eq:var}) with a link ``$X_{t-\tau}\to Y_t$'' and $\tau > 0$ can be written as
\begin{align} 
X_t &= \sum_{i=1}^{N_X} a_i Z^i_{t-h_i} + \varepsilon_{X,t}  \label{subX}  \\
Y_t &= c X_{t-\tau} + \sum_{i=1}^{N_Y} b_i W^i_{t-g_i} + \varepsilon_{Y,t} \label{subY}
\end{align}
with parents
\begin{align}
Z^i_{t-h_i} &\in \mathcal{P}_{X_{t}}~~~\text{for}~~~i =1,\ldots,N_X,\\
W^i_{t-g_i} &\in \mathcal{P}_{Y_{t}}{\setminus}\{X_{t-\tau}\}~~~\text{for}~~~i =1,\ldots,N_Y.
\end{align}
Here the coefficient $c$ corresponds to the entry $\Phi(\tau)_{YX}$.

To simplify notation, Eqns.~(\ref{subX},\,\ref{subY}) are expressed in vector notation
\begin{align}  \label{eq:var_vector}
X_t &= \mathbf{Z}_t A + \varepsilon_{X,t}  \nonumber  \\
Y_t &= cX_{t-\tau} + \mathbf{W}_t B + \varepsilon_{Y,t}  
\end{align}
where $X_t,Y_t$ are scalar random processes, $A$ and $B$ are the coefficient vectors, and $\mathbf{Z}_t,\, \mathbf{W}_t$ are possibly multivariate random processes of dimension $ N_X$ and  $N_Y$ respectively,
\begin{align}
\mathbf{Z}_t&=(Z^1_{t-h_1},\ldots,Z^{N_X}_{t-h_{N_X}}),  \\
\mathbf{W}_t &= (W^1_{t-g_1},\ldots,Z^{N_Y}_{t-g_{N_Y}}).
\end{align}
In the following, $t$ and $\tau$ will be dropped for ease of notation.

For the cross correlation between $X$ and $Y$ at lag $\tau$, the covariance $E[Y^\top X]$ and the variances $E[Y^\top Y]$ and $E[X^\top X]$ are needed.
While in the covariance expression Eq.~(\ref{eq:anacov}) the dependencies are rather hidden, the vector notation allows to derive them simply by directly plugging in Eqns.~(\ref{eq:var_vector}) into the covariances and using only $E[\mathbf{W}^\top \varepsilon_Y]=E[\mathbf{Z}^\top \varepsilon_Y]=E[\mathbf{Z}^\top \varepsilon_X]=0$ since $\varepsilon_Y$ and $\varepsilon_X$ are i.i.d. processes independent from the past parents. Then the (co-)variances can be written in a compact way:
\begin{align}
E[Y^\top X] &= c \sigma_X^2  \nonumber\\
                &~~~+ c A^\top E[\mathbf{Z}^\top\mathbf{Z}]A \nonumber \\
                &~~~+ B^\top E[ \mathbf{W}^\top \mathbf{Z}] A \nonumber\\
               &~~~+ B^\top E[\mathbf{W}^\top \varepsilon_X],\\
E[Y^\top Y] &= \sigma_Y^2 + c^2 \sigma_X^2  \nonumber\\
    &~~~+ c^2 A^\top E[\mathbf{Z}^\top\mathbf{Z}]A + B^\top E[ \mathbf{W}^\top \mathbf{W}] B  \nonumber\\
    &~~~+ c\left( B^\top E[ \mathbf{W}^\top \mathbf{Z}] A + A^\top E[ \mathbf{Z}^\top \mathbf{W}] B \right)  \nonumber \\
    &~~~+ c \left(        B^\top E[\mathbf{W}^\top\varepsilon_X] + E[\varepsilon_X^\top\mathbf{W}] B \right), \\
E[X^\top X] &= \sigma_X^2 \nonumber\\
               &~~~+ A^\top E[\mathbf{Z}^\top\mathbf{Z}]A.
\end{align}
One can see, that the covariance $E[Y^\top X]$ not only depends on the coefficient $c$, but also on the variance of the parents $\mathbf{Z}$ of $X$, the covariance among the parents of $X$ and $Y$ and the covariance of the innovation $\varepsilon_X$ with the parents $\mathbf{W}$ of $Y$.

Also in this interpretation, we find that the value of the cross correlation cannot easily be related to the coefficient $c$ of the link between $X$ and $Y$ in the time series graph and depends on the multiple interactions between the parents of $X$ and $Y$ in the multivariate process.

\section{Partial correlation measure MIT of multivariate autoregressive process} \label{par_corr}
\subsection{Definitions}
The knowledge of the (linear) conditional independence structure of the data encoded in the time series graph can be used to define a certain partial correlation measure with a straightforward graph-theoretical interpretation. 

Partial correlation can be defined in the framework of regression analysis.
If one regresses two variables $X,\,Y$ on the same regressors $\mathbf{U}$,  the cross correlation between the residuals
\begin{align} \label{eq:residuals}
X_{\mathbf{U}} &\equiv X - \mathbf{U} (E[\mathbf{U}^\top \mathbf{U}])^{-1} E[\mathbf{U}^\top X] \nonumber\\
Y_{\mathbf{U}} &\equiv Y - \mathbf{U} \underbrace{(E[\mathbf{U}^\top \mathbf{U}])^{-1} E[\mathbf{U}^\top Y]}_{\substack{\text{regression coefficient}\\\text{(vector) } \mathbf{R}}} .
\end{align}
is the partial correlation 
\begin{align} \label{eq:def_parcorr}
\rho(YX|\mathbf{U}) &= \frac{E[Y_{\mathbf{U}}^\top X_{\mathbf{U}}]}{\sqrt{E[Y_{\mathbf{U}}^\top Y_{\mathbf{U}}]}\sqrt{E[X_{\mathbf{U}}^\top X_{\mathbf{U}}}}.
\end{align}
Note, that this measure is not to be confused with the partial autocorrelation \citep{brockwell2009time, von2002statistical}.

The partial correlation measure introduced now is based on the parents $Y$ \emph{and} the parents of $X$. 
\begin{mydef}
For two components $X,\,Y$ of a stationary multivariate discrete-time process $\mathbf{X}$ with parents $\mathcal{P}_{Y_t}$ and $\mathcal{P}_{X_t}$ in the associated time series graph and $\tau>0$,
\begin{align} \label{eq:def_mit}
 \rho^{\rm MIT}_{X{\to}Y}(\tau) &\equiv \rho(X_{t-\tau};Y_t| \mathcal{P}_{Y_t} {\setminus}\{X_{t-\tau}\},\mathcal{P}_{X_{t-\tau}} ).
\end{align}
\end{mydef}
The name MIT, short for \emph{momentary information transfer}, is used in analogy to the general case described in \cite{Runge2012b}, which in the linear case should be understood as \emph{momentary variance transfer}. The attribute \textit{momentary} \citep{Pompe2011} is used because MIT measures the variance of the ``moment'' $t-\tau$ in $X$ that is transferred to $Y_t$. $\rho^{\rm MIT}$ quantifies how much the variability in $X$ at the exact lag $\tau$ \emph{directly} influences $Y_t$, irrespective of the pasts of $X_{t-\tau}$ and $Y_t$. 
One can also define a contemporaneous MIT, which in the linear case is equivalent to the inverse covariance of the residuals after regressing each process on its parents \citep{Runge2012b}.


\subsection{Properties: Linear coupling strength autonomy theorem}
As in Sect.~\ref{sect:cc_deps_par} for the cross correlation, we now derive the dependencies of the partial correlation MIT on the coefficients of a vector autoregressive model Eq.~(\ref{eq:var_vector}). The equations for the subprocess $Y$ can be written as
\begin{align}  \label{eq:var_vector_reg}
Y_t &= \mathbf{W}_t B + \varepsilon_{Y,t},
\end{align}
where $X$ and the coefficient $c$ occurring in Eq.~(\ref{eq:var_vector}) is collapsed into $\mathbf{W}$ and $B$, respectively.  
\begin{mylemma}
For the autoregressive model Eq.~(\ref{eq:var_vector_reg}), a multivariate regression for the dependent variable $Y$ on $\mathbf{U}=(\mathbf{W},\,\mathbf{V})$, where $\mathbf{V}$ are other regressors that are not part of the parents, i.e., $\mathbf{V}\cap \mathbf{W}=\emptyset$ gives
\begin{align}
\left(\begin{array}{c} \mathbf{R}_W \\ \mathbf{R}_V \end{array} \right) &= \left(\begin{array}{c} \mathbf{B} \\ 0 \end{array} \right).
\end{align}
\end{mylemma}
\noindent
The proof is given in the appendix. For the partial correlation MIT, the dependencies are slightly more complex.
\begin{mythm} \label{thm}
For the autoregressive model Eq.~(\ref{eq:var}), written in vector notation as Eq.~(\ref{eq:var_vector}), the partial correlation $\rho^{\rm MIT}_{X{\to}Y}(\tau)$ given by Eq.~(\ref{eq:def_mit}) written in vector notation as in Eq.~(\ref{eq:def_parcorr}) with $\mathbf{U}=(\mathbf{W},\,\mathbf{Z})$ is comprised of the covariances and variances
\begin{align} \label{eq:mit_theorem}
E[Y_{\mathbf{U}} ^\top X_{\mathbf{U}} ] &= c \sigma_X^2  \nonumber\\
               &~~~- c E[\varepsilon_X^\top\mathbf{W}] S_Z^{-1}   E[\mathbf{W}^\top \varepsilon_X] \nonumber\\
E[Y_{\mathbf{U}} ^\top Y_{\mathbf{U}} ] &= \sigma_Y^2 + c^2 \sigma_X^2  \nonumber\\
    &~~~- c^2 E[\varepsilon_X^\top\mathbf{W}] S_Z^{-1}   E[\mathbf{W}^\top \varepsilon_X]\nonumber\\
E[X_{\mathbf{U}} ^\top X_{\mathbf{U}} ] &= \sigma_X^2  \nonumber\\
               &~~~- E[\varepsilon_X^\top\mathbf{W}] S_Z^{-1}   E[\mathbf{W}^\top \varepsilon_X].
\end{align}
where $S_Z$ denotes the Schur complement
\begin{align}
S_Z &=  E[\mathbf{W}^\top \mathbf{W}] - E[\mathbf{W}^\top \mathbf{Z}] (E[\mathbf{Z}^\top\mathbf{Z}])^{-1} E[\mathbf{Z}^\top \mathbf{W}].
\end{align}
\end{mythm}
The proof is given in the appendix. The (co-)variances are comprised of two parts. The first one is simply the cross correlation between $\varepsilon_{X,t-\tau}$ and $\varepsilon_{Y,t}+c\,\varepsilon_{X,t-\tau}$. The second part is due to dependencies between $\varepsilon_{X,t-\tau}$ and the parents of $Y_t$ and non-zero only under certain conditions.

More precisely, the Schur complement $S_Z$ can be interpreted as the conditional variance of $\mathbf{W}$ given $\mathbf{Z}$. On the other hand, the covariance $E[\varepsilon_X^\top \mathbf{W}]$ can best be interpreted in the framework of time series graphs. In terms of the coefficient path matrices $\Psi$ and the innovation's covariance $\Sigma$ it can be written as: 
\begin{align}
(E[\varepsilon_X^\top \mathbf{W}])_i = \sum^N_{r=1} \Psi_{W_ir}(\tau-g_i) \Sigma_{rX}.
\end{align}
This relation is derived in the appendix. $(E[\varepsilon_X^\top \mathbf{W}])_i$ is the linear combination of all paths of length $\tau-g_i$ emanating from $X_t$ or $\mathbf{X}_{r,t}$ with $\Sigma_{rX} \neq 0$ to $W^i_{t+\tau-g_i}$. 
It will be shown, that it can be understood as a ``sidepath'' covariance and is zero if there are no such paths.
Then, for $E[\varepsilon_X^\top \mathbf{W}]=0$, the $\rho^{\rm MIT}$ becomes
\begin{align} \label{eq:mit_nosidepaths}
\rho_{X{\to}Y}^{\rm MIT} &= \frac{c \sigma_X}{\sqrt{ \sigma_Y^2 + c^2 \sigma_X^2}}.
\end{align}
Thus, if there are no sidepaths, the partial correlation measure MIT of a link ``$X_{t-\tau}\to Y_t$'' solely depends on the coefficient matrix entry $\Phi_{YX}(\tau)$ and the innovation's variances $\sigma^2_X$ and $\sigma^2_Y$. The MIT of an autoregressive process is, therefore, much better interpretable than the cross correlation as analysed in Sects.~\ref{sect:cc_deps_paths} and \ref{sect:cc_deps_par} since its value is attributable to the interaction between $\mathbf{X}^j$ and $\mathbf{X}^i$ alone, i.e., the link ``$\mathbf{X}^j_{t-\tau}\to\mathbf{X}^i_t$'' in the time series graph of $\mathbf{X}$.
This theorem is the linear version of the \emph{coupling strength autonomy theorem} that treats the general nonlinear case in the information-theoretic framework \citep{Runge2012b}.


\subsection{Alternative measures}
The graph-theoretic perspective invites to define related measures that capture different aspects of the dependency between two components in a multivariate process.

For example, we can also choose either one of the parents as a condition, which -- dropping the attribute ``momentary'' -- leads to the \emph{information transfers} ITY and ITX
\begin{align} 
\rho_{X \to Y}^{\rm ITY}(\tau) &\equiv \rho(X_{t-\tau};Y_t|\mathcal{P}_{Y_t}{\setminus}\{X_{t-\tau}\}) \label{eq:def_py}, \\
\rho_{X \to Y}^{\rm ITX}(\tau) &\equiv \rho(X_{t-\tau};Y_t|\mathcal{P}_{X_{t-\tau}})  \label{eq:def_px}.
\end{align}
ITY only conditions out the influence of the parents of $Y$, but includes the aggregated influence of the parents of $X$. Like MIT it is non-zero only for (Granger-) causal dependent nodes and used in the algorithm to estimate the time series graph \citep{Runge2012prl,Runge2013}. ITX, on the other hand, measures the part of variance originating in $X_{t{-}\tau}$ that reaches $Y_t$ on any path and is, thus, not a `causal' measure of direct dependence, yet in many situations we might only be interested in the effect of $X$ on $Y$, no matter how this influence is mediated. 

For the case of sidepaths with $E[\varepsilon_X^\top \mathbf{W}]\neq 0$ the (co-)variances in Eq.~(\ref{eq:mit_theorem}) depend on an additional term.
As an example where one parent $W_k$ of $Y$ (apart from $X$) depends on $X$, consider the following model:
\begin{align} \label{sidepath_model}
X_t &= \sum_{i=1}^{N_X} a_i Z^i_{t-h_i} + \varepsilon_{X,t}   \\
W^k_t &= d X_{t-1} + \varepsilon_{W_k,t}    \\
Y_t &= c X_{t-2} + b_k W^k_{t-1} + \nonumber\\
&~~~+ \sum_{i=1\neq k}^{N_Y} b_i W^i_{t-g_i} + \varepsilon_{Y,t}  
\end{align}
where for all $i\neq k:~(E[\varepsilon_X^\top \mathbf{W}])_i=0$ and also assume that additionally for all $i \neq X:~\Sigma_{iX}=0$. 
As derived in the appendix, MIT is then
\begin{align}
\rho^{\rm MIT}_{X{\to}Y}(\tau)&= \rho(X_{t-\tau};Y_t~|~\mathbf{W}_t,\mathbf{Z}_{t-\tau}) \nonumber\\
& =\frac{c \sigma_{W_k}^2 \sigma_{X}^2}{\sqrt{c^2 \sigma_{W_k}^2 \sigma_{X}^2+\left(\sigma_{W_k}^2+d^2 \sigma_{X}^2\right) \sigma_{Y}^2}}.
\end{align}
Thus, the MIT depends not only on $c$, but also on all the coefficients along the paths $\Psi_{W_kX}(\tau-g_k)$, here only $d$, and on the residual variance of $W_k$ given $\mathbf{Z}$.

This example points to the suggestion, that it might be more appropriate to ``leave open'' \emph{all} paths from $X_{t-\tau}$ to $Y_t$ by excluding from the conditions those parents of $Y_t$ that are depending on $X_{t-\tau}$. Then the possible paths of variance transfer are either via the direct link ``$X_{t-\tau}\to~Y_t$'' or via the sidepaths ``$X_{t-\tau}\mathop{}_{\--}^{\to}\cdots \to \cdots \to Y_t$'' (the symbol ``$\mathop{}_{\--}^{\to}$'' denotes that the sidepath can start from $X_{t-\tau}$ either directed or contemporaneous, while the subsequent links of the path can only be directed). To isolate all of these paths, we suggest to additionally condition on the parents of the intermediate nodes on these sidepaths. These nodes can be characterised by 
\begin{align}
\mathcal{A}^\star_{Y_t}\equiv \{& W^k_{t-\tau_k} \in \mathcal{A}_{Y_t}{\setminus} \{X_{t-\tau},\mathcal{P}_{X_{t-\tau}}\} :\nonumber\\
                         &\rho(W^k_{t-\tau_k}; X_{t-\tau}|\mathcal{P}_{X_{t-\tau}})\neq 0 \},
\end{align}
where $\mathcal{A}_{Y_t}$ denotes the ancestors of $Y_t$, i.e., the set of nodes with a directed path towards $Y_t$ \citep{Eichler2011}.  We call the modified MIT MITS, where ``S'' stands for ``sidepath,''
\begin{align} \label{eq:mit_star}
\rho^{\rm MITS}_{X{\to}Y}(\tau)&\equiv \rho(X_{t-\tau};Y_t~|~ \{ \mathcal{P}_{Y_t},\mathcal{P}(\mathcal{A}^\star_{Y_t})\}{\setminus} \{ \mathcal{A}^\star_{Y_t},X_{t-\tau}\},\mathcal{P}_{X_{t-\tau}} ).
\end{align}
In our sidepath example Eq.~(\ref{sidepath_model}) for the simpler special case $a_i=0~\forall i$ and $b_i=0~\forall~i \neq k$, MITS evaluates to \citep{Runge2012b}
\begin{align}
\rho^{\rm MITS}_{X{\to}Y}(\tau) &=\frac{(c+ d b_k) \sigma_X}{\sqrt{c^2 \sigma_X^2+b_k \left(b_k \sigma_{W_k}^2+d (2 c+b_k d) \sigma_X^2\right)+\sigma_Y^2}}.
\end{align}
Here the factor $c+d b_k$ is the covariance along both paths, which can also vanish for $c=-d b_k$, and seems like a more appropriate representation of the coupling between $X_{t-2}$ and $Y_t$.


\section{Analysis of sampling distributions} \label{numerics}
\begin{figure}[!t]
\begin{center}
\includegraphics[width=.8\columnwidth]{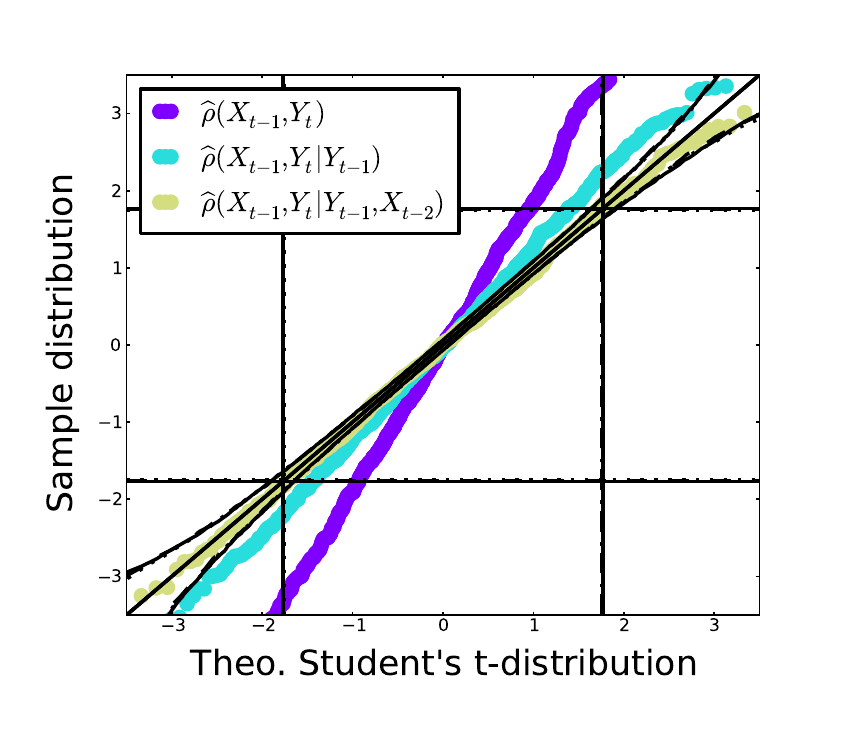}
\end{center}
\caption[]{Quantile (``$q-q$''-) plots of sample estimates of cross correlation and the partial correlation measures ITY and MIT plotted against the Student's t-distribution for different degrees of freedom (dotted line for $q=0$, dashed line for $q=1$ and solid line for $q=2$; the lines are almost identical). The diagonal line (with 90\% confidence intervals) indicates a perfect match of theoretical and empirical distributions and the horizontal and vertical black lines denote the 5\% and 95\% quantiles of the theoretical distributions for different degrees of freedom.}
\label{fig:qqplot}
\end{figure}



In this section we study the properties of the sample estimate $\widehat{\rho_{X{\to}Y}^{\rm MIT}}$ of the MIT partial correlation. 
It is known, that the distribution of the partial correlation coefficient is the same as that of the cross correlation coefficient with the degrees of freedom reduced by the cardinality of the set of conditions $q$ \citep{Fisher1924}. Therefore, the distribution of
\begin{align}
\widehat{t}(YX|\mathbf{U}) &= \widehat{\rho}(YX|\mathbf{U}) \sqrt{\frac{n-2-q}{1-\widehat{\rho}(YX|\mathbf{U})^2}}
\end{align}
is Student's-$t$ with $n-2-q$ degrees of freedom with $q$ being the dimension of $\mathbf{U}$. In the case of MIT $q=|\{\mathcal{P}_{Y_t} {\setminus}\{X_{t-\tau}\},\mathcal{P}_{X_{t-\tau}}\}|$. 

The assumptions underlying this result are Gaussianity and, importantly, independent and identically distributed samples. This assumption is, however, violated in many practical cases, especially for serially dependent samples with non-zero autocorrelations. Consider the model Eq.~(\ref{eq:ar_model_matrix}) for $c=0$ and strong autocorrelations $a=b=0.9$ and where we assume the innovations to be uncorrelated, i.e., $\Sigma$ is diagonal. The two processes are, therefore, independent, but the samples are serially dependent. As shown in Fig.~\ref{fig:qqplot} for the cross correlation this effectively reduces the degrees of freedom \citep{Chatfield2003} and leads to an ``inflated'' sampling distribution. 

Since Theorem \ref{thm} implies, that MIT ``filters out'' also autocorrelation, we expect that, conversely to the sample estimate of the cross correlation, the MIT estimator is not ``inflated'' by autocorrelation. More precisely, since the condition on the parents removes the dependency of $X$ and $Y$ on the past samples, the residuals $X_{\mathbf{U}}$ and $Y_{\mathbf{U}}$ given by Eq.~(\ref{eq:residuals}) for a regression on both parents $\mathbf{U}=\{\mathcal{P}_{Y_t} {\setminus}\{X_{t-\tau}\},\mathcal{P}_{X_{t-\tau}}\}$ are
\begin{align}
X_{\mathbf{U},t} &= \varepsilon_{X,t}\\
Y_{\mathbf{U},t} &= c\varepsilon_{X,t} + \varepsilon_{Y,t}
\end{align}
 and therefore indeed serially independent since both $\varepsilon_{X,t}$ and $\varepsilon_{Y,t}$ are independent in time. Note, that this only holds for links ``$X_{t-\tau}\to Y_t$'' without sidepaths as discussed in the previous section.
We also test the distribution of the partial correlation ITY defined in Eq.~(\ref{eq:def_py}) where only the parents of $Y$ are conditioned out. Here the residuals are \emph{not} independent and we expect the distribution to be still broadened due to less effective degrees of freedom. For model Eq.~(\ref{eq:ar_model_matrix}) the parents are $\mathcal{P}_{Y_t} \setminus\{X_{t-\tau}\} = Y_{t-1}$ and $\mathcal{P}_{X_{t-\tau}} = X_{t-2}$ and $\tau=1$.

Figure~\ref{fig:qqplot} shows the quantile plots of the empirical distributions simulated with time series length $T=20$ plotted against the Student's t-distribution with $q=0$ for $\widehat{\rho}$, $q=1$ for $\widehat{\rho_{X{\to}Y}^{\rm ITY}}$ and $q=2$ for $\widehat{\rho_{X{\to}Y}^{\rm MIT}}$. 
The plots demonstrate, that the cross correlation is strongly ``inflated'', ITY is still affected and only MIT can be well described by the theoretical distribution within the confidence bounds, independent of the strength of autocorrelation. 

This feature can be used for independence tests since it allows for a more accurate significance test. 
Note, however, that first the time series graph has to be estimated to infer the parents to condition on. In \cite{Runge2013} the measure ITY is used in the estimation of the time series graph and we suggest to subsequently test the inferred links with MIT to fully account for autocorrelations and dependencies also from parents of $X$.

\section{Application to climatological time series} \label{application}
As a climatological application we study two indices of monthly sea surface temperature anomalies \citep{rayner2003global} for the period 1950 -- 2012. NINO is the time series of the spatial average over the Nino34 region (5N-5S and 170-120W) in the East Pacific and TNA is the tropical North Atlantic index \citep{Enfield1999} averaged over (5.5--23.5N and 15--57.5W). 

\begin{figure}[!t]
\begin{center}
\includegraphics[width=\columnwidth]{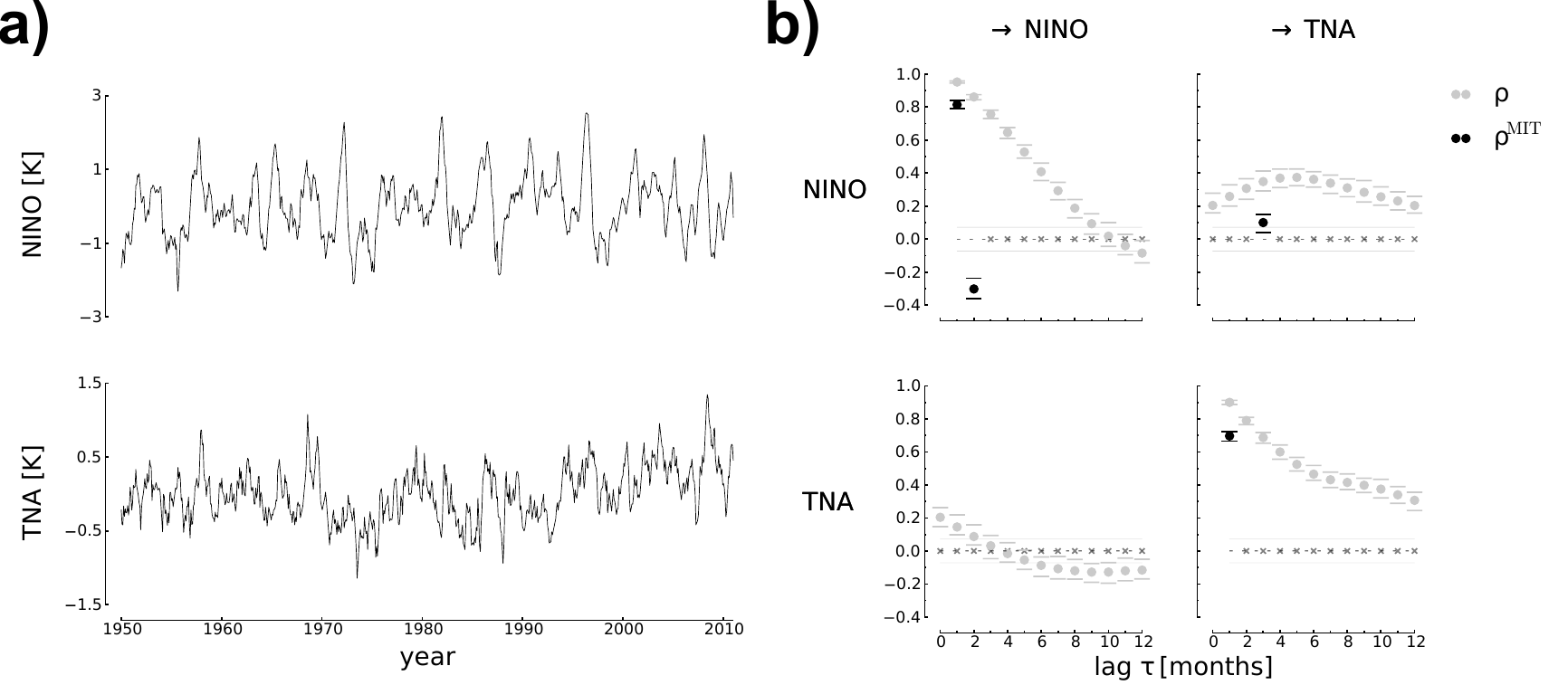}
\end{center}
\caption[]{
(a) Time series and (b) correlations and partial correlations of climatological example. The matrix of lag functions shows the (auto-)correlations (light gray) and the value of MIT (black), where non-significant links are marked by gray crosses. The horizontal gray line denotes the two-sided 95\%-significance level for the (auto-)correlations. The errorbars mark the 90\% confidence interval estimated from a bootstrap test. For example, the upper right plot shows the lagged cross correlation function $\rho(\text{NINO}_{t-\tau};\text{TNA}_t)$ for $\tau\geq 0$ in light gray and the MIT value at the significant link ``$\text{NINO}_{t-3}\to\text{TNA}_t$'' in black. Note, that for autocorrelations (on the diagonal) the zero-lag is not drawn. 
}
\label{fig:climate}
\end{figure}

Figure~\ref{fig:climate} shows the time series and (partial) correlations. The time series graph was estimated using the PC-algorithm \citep{Spirtes2000} as described in \cite{Runge2012prl,Runge2013} with the theoretical significance test discussed above at the (two-sided) level $\alpha=95\%$. The estimated time series graph is comprised of a coupling link ``$\text{NINO}_{t-3}\to\text{TNA}_t$'' and autodependency links at lag 1 and 2 in NINO and only at lag 1 in TNA. On the other hand, the auto- and cross correlation lag functions shown in gray feature significant links for a large range of lags with a maximum of the cross correlation lag function $\rho(\text{NINO}_{t-\tau};\text{TNA}_t)$ at lag $\tau=5$. This shift of the lag function's maximum is further investigated in \cite{Runge2013}. 
Also the cross correlation value $\rho=0.35\pm0.05$ at lag $\tau=3$ is significantly larger than $\rho^{\rm MIT}(\tau=3)=0.10\pm0.05$ (the ``$\pm$'' values correspond to the 90\% confidence interval estimated from a bootstrap test \citep{Runge2013}). 

The strong autodependency links with MIT values of $(0.8,\,-0.3)$ for lags 1 and 2 in NINO and 0.7  for lag 1 in TNA explain these `significant' cross correlation values at most lags, which according to Eq.~(\ref{eq:anacov}), are due to the common driver effect of past nodes (Fig.~\ref{fig:causality}(b)) or the indirect causal effect due to intermediate lags (Fig.~\ref{fig:causality}(a)). 
On the other hand, since there are no sidepaths here, the small MIT value reflects only the contributions from the coupling link and the residual's variances according to Eq.~(\ref{eq:mit_nosidepaths}). The small value of MIT shows, that the actual coupling mechanism by which NINO influences TNA is quite weak, but due to strong autocorrelations the overall contribution to TNA's variance is larger becoming maximal in the peak at lag 5. In \cite{Runge2013} this Pacific -- Atlantic interaction is climatologically discussed.

\section{Conclusions}
With the goal to investigate how the value of cross correlation can be interpreted, we analysed how the cross correlation between two components of a multivariate autoregressive model depends on the model's coefficients. These dependencies can well be understood within the framework of time series graphs showing that the value of cross correlation at a certain lag stems from a superposition of paths from past and intermediate nodes in the graph. These complex dependencies on the model's coefficients make it hard to interpret the cross correlation as a measure of the strength of association between the two components alone. 

On the other hand, for the recently introduced partial correlation measure MIT we prove a simple formula depending solely on the coefficient belonging to the coupling lag of the two variables and the variance of their innovations. MIT, thus, allows to separate the effect of other links, like strong autocorrelations, from the actual coupling link, making it better interpretable than cross correlation. We also suggest related measures that capture different aspects of the dependency between two components in a well interpretable way. Additionally, an analysis of the sample estimate of MIT shows, that it is not `inflated' by autocorrelations like cross correlation and, thus, suitable for significance tests that assume temporally independent samples.

On {\tt http://tocsy.pik-potsdam.de/tigramite.php} we provide a Python program with a graphical user interface to estimate the time series graph and the partial correlation measures ITY and MIT as well as their information-theoretic counterparts.

\section*{\small Acknowledgment}
We appreciate the support by the German National Academic Foundation (Studienstiftung) and DFG grant No. KU34-1 and thank Jobst Heitzig for helpful comments on an earlier version of the manuscript.

\appendix
\section{Appendix}
\subsection{Proof that parents as regressors yields model coefficients}
For the model Eq.~(\ref{eq:var}) any regression of $Y$ on regressors that include the parents $\mathcal{P}_{Y_t}$ yields the corresponding coefficients in $\Phi$ for the parents and zeros for non-parents. More precisely, first the dependencies of a subprocess $Y\in \mathbf{X}$ can be written as
\begin{align} 
Y_t &= \sum_{i=1}^{N_Y} b_i W^i_{t-g_i} + \varepsilon_{Y,t} \label{subYreg}
\end{align}
with parents
\begin{align}
W^i_{t-g_i} &\in \mathcal{P}_{Y_{t}}.
\end{align}
To simplify notation, Eq.~(\ref{subYreg}) is expressed in vector notation
\begin{align}
Y_t &= \mathbf{W}_t \mathbf{B} + \varepsilon_{Y,t}  \label{Yreg},
\end{align}
where $\mathbf{B}=(b_1,\ldots)$ is the coefficient vector and $\mathbf{W}_t$ is a possibly multivariate random process of dimension $N_Y$, on which $Y$ depends at lags $g_1,\ldots,g_{N_Y}$,
\begin{align}
\mathbf{W}_t &= (W^1_{t-g_1},\ldots,W^{N_Y}_{t-g_{N_Y}}).
\end{align}
In the following, $t$ and $\tau$ will be dropped for ease of notation.

Then a regression on $\mathbf{U}=(\mathbf{W},\,\mathbf{V})$, where $\mathbf{V}$ are other regressors that are not part of the parents, i.e., $\mathbf{V}\cap \mathbf{W}=\emptyset$ gives the coefficient vector
\begin{align} \label{eq:reg_inv}
\left(\begin{array}{c} \mathbf{R}_W \\ \mathbf{R}_V \end{array} \right) &\equiv (E[\mathbf{U}^\top \mathbf{U}])^{-1} E[\mathbf{U}^\top Y] \nonumber\\
   &= \left( \begin{array}{cc} E[\mathbf{W}^\top\mathbf{W}] & E[\mathbf{W}^\top\mathbf{V}] \\ E[\mathbf{V}^\top\mathbf{W}] & E[\mathbf{V}^\top\mathbf{V}] \end{array} \right)^{-1} \left(\begin{array}{c} E[\mathbf{W}^\top Y] \\ E[\mathbf{V}^\top Y] \end{array} \right).
\end{align}
Now one can prove that 
\begin{align}
\left(\begin{array}{c} \mathbf{R}_W \\ \mathbf{R}_V \end{array} \right) &= \left(\begin{array}{c} \mathbf{B} \\ 0 \end{array} \right),
\end{align}
which implies that any multivariate regression which contains the parents as regressors will recover the coefficients of the underyling model. 

To prove this relation, the inverse can be treated via the matrix inversion lemma
\begin{align} \label{eq:inversion_lemma}
&\left( \begin{array}{cc} E[\mathbf{W}^\top\mathbf{W}] & E[\mathbf{W}^\top\mathbf{V}] \\ E[\mathbf{V}^\top\mathbf{W}] & E[\mathbf{V}^\top\mathbf{V}] \end{array} \right)^{-1} \nonumber\\
&=\left( \begin{smallmatrix}
S^{-1}_V & - (E[\mathbf{W}^\top \mathbf{W}])^{-1} E[\mathbf{W}^\top \mathbf{V}] S^{-1}_W \\
- (E[\mathbf{V}^\top \mathbf{V}])^{-1} E[\mathbf{V}^\top \mathbf{W}] S^{-1}_V & S_W^{-1}  \end{smallmatrix} \right),
\end{align}
where $S_{\cdot}$ denotes the Schur complements
\begin{align}
S_V &=  E[\mathbf{W}^\top \mathbf{W}] - E[\mathbf{W}^\top \mathbf{V}] (E[\mathbf{V}^\top\mathbf{V}])^{-1} E[\mathbf{V}^\top \mathbf{W}] \\
S_W &=  E[\mathbf{V}^\top \mathbf{V}] - E[\mathbf{V}^\top \mathbf{W}] (E[\mathbf{W}^\top\mathbf{W}])^{-1} E[\mathbf{W}^\top \mathbf{V}].
\end{align}
$S_V$ can be interpreted as the conditional variance of $\mathbf{W}$ given $\mathbf{V}$. $S^{-1}_V$ can be further transformed using the Woodbury matrix identity
\begin{align} \label{eq:woodbury}
S^{-1}_V &= (E[\mathbf{W}^\top \mathbf{W}])^{-1} - (E[\mathbf{W}^\top \mathbf{W}])^{-1} E[\mathbf{W}^\top \mathbf{V}]\times \nonumber\\
&\times \underbrace{(-E[\mathbf{V}^\top \mathbf{V}] + E[\mathbf{V}^\top \mathbf{W}] (E[\mathbf{W}^\top\mathbf{W}])^{-1} E[\mathbf{W}^\top \mathbf{V}])^{-1}}_{=(-S^{-1}_W)}   \times\nonumber \\
&~~~~~~\times E[\mathbf{V}^\top \mathbf{W}] (E[\mathbf{W}^\top \mathbf{W}])^{-1}.
\end{align}
The covariance vector in Eq.~(\ref{eq:reg_inv}) can be simplified by
\begin{align}
\left(\begin{array}{c} E[\mathbf{W}^\top Y] \\ E[\mathbf{V}^\top Y] \end{array} \right) &= \left(\begin{array}{c} E[\mathbf{W}^\top \mathbf{W}] B + E[\mathbf{W}^\top \varepsilon_Y] \\ E[\mathbf{V}^\top \mathbf{W}]B +  E[\mathbf{V}^\top \varepsilon_Y] \end{array} \right) \\
&= \left(\begin{array}{c} E[\mathbf{W}^\top \mathbf{W}] B \\ E[\mathbf{V}^\top \mathbf{W}]B \end{array} \right),
\end{align}
where $ E[\mathbf{V}^\top \varepsilon_Y]=E[\mathbf{W}^\top \varepsilon_Y]=0$ because $\varepsilon_Y$ is independent of past processes.
Then the regression coefficient $\mathbf{R}_W$ given by
\begin{align}
\mathbf{R}_W &= S_V^{-1}  E[\mathbf{W}^\top \mathbf{W}] \mathbf{B} - \nonumber\\
&~~~~- ( E[\mathbf{W}^\top \mathbf{W}] )^{-1}  E[\mathbf{W}^\top \mathbf{V}] S_W^{-1}  E[\mathbf{V}^\top \mathbf{W}] \mathbf{B}  
\end{align}
can be simplified by inserting Eq.~(\ref{eq:woodbury}) from which it follows that
\begin{align}
&S_V^{-1}  E[\mathbf{W}^\top \mathbf{W}] \mathbf{B} =  \nonumber\\
&\mathbf{B} + (E[\mathbf{W}^\top \mathbf{W}])^{-1} E[\mathbf{W}^\top \mathbf{V}] S_W^{-1} E[\mathbf{V}^\top \mathbf{W}] \mathbf{B},
\end{align}
and thus $\mathbf{R}_W=\mathbf{B}$ which proves the first part of the claim.

To prove the second part, now the analogue of Eq.~(\ref{eq:woodbury}) for $S^{-1}_W$ is inserted into
\begin{align}
\mathbf{R}_V &= S_W^{-1}  E[\mathbf{V}^\top \mathbf{W}] \mathbf{B} - \nonumber\\
&~~~~- ( E[\mathbf{V}^\top \mathbf{V}] )^{-1}  E[\mathbf{V}^\top \mathbf{W}] S_V^{-1}  E[\mathbf{W}^\top \mathbf{W}] \mathbf{B},
\end{align}
from which using
\begin{align}
&S_W^{-1}  E[\mathbf{V}^\top \mathbf{W}] \mathbf{B} = (E[\mathbf{V}^\top \mathbf{V}])^{-1} E[\mathbf{V}^\top \mathbf{W}] \mathbf{B}\times \nonumber\\
&~~~~~~~\times (E[\mathbf{V}^\top \mathbf{V}]])^{-1} E[\mathbf{V}^\top \mathbf{W}] S_V^{-1} \times \nonumber \\
&~~~~~~~\times E[\mathbf{W}^\top \mathbf{V}] (E[\mathbf{V}^\top \mathbf{V}])^{-1}E[\mathbf{V}^\top \mathbf{W}] \mathbf{B},
\end{align}
and
\begin{align}
E[\mathbf{W}^\top \mathbf{V}] (E[\mathbf{V}^\top \mathbf{V}])^{-1}E[\mathbf{V}^\top \mathbf{W}] = E[\mathbf{W}^\top \mathbf{W}] - S_V
\end{align}
one arrives at $\mathbf{R}_V=0$.

\subsection{Proof of linear coupling strength autonomy theorem}
First $X$ and $Y$ are regressed on $\mathbf{U} =(\mathbf{W},\,\mathbf{Z})$ yielding the residuals
\begin{align}
Y_{\mathbf{U}} &\equiv Y - \mathbf{U} (E[\mathbf{U}^\top \mathbf{U}])^{-1} E[\mathbf{U}^\top Y] \\
X_{\mathbf{U}} &\equiv X - \mathbf{U} (E[\mathbf{U}^\top \mathbf{U}])^{-1} E[\mathbf{U}^\top X] .
\end{align}
Then the covariance and variances are
\begin{align}
E[Y_{\mathbf{U}} ^\top X_{\mathbf{U}} ] &= E[Y^\top X]-    \nonumber\\
               &~~~- E[Y^\top \mathbf{U}] (E[\mathbf{U}^\top \mathbf{U}])^{-1} E[\mathbf{U}^\top X]\\
E[Y_{\mathbf{U}} ^\top Y_{\mathbf{U}} ] &= E[Y^\top Y] -   \nonumber\\
               &~~~- E[Y^\top \mathbf{U}] (E[\mathbf{U}^\top \mathbf{U}])^{-1} E[\mathbf{U}^\top Y]\\
E[X_{\mathbf{U}} ^\top X_{\mathbf{U}} ] &= E[X^\top X] -   \nonumber\\
               &~~~- E[X^\top \mathbf{U}] (E[\mathbf{U}^\top \mathbf{U}])^{-1} E[\mathbf{U}^\top X].
\end{align}
The covariance can be evaluated as follows. First, writing
\begin{align}
X &= \mathbf{U} \left( \begin{smallmatrix} 0\\ A \end{smallmatrix} \right) + \varepsilon_X \\
Y & = \mathbf{U} \left( \begin{smallmatrix} \mathbf{B}\\ c \mathbf{A} \end{smallmatrix} \right) +c \varepsilon_X + \varepsilon_Y
\end{align}
the covariance $E[Y^\top X]$ is expressed in terms of $\mathbf{U}$ as
\begin{align}
&E[Y^\top X] = \nonumber\\
& \left( \mathbf{B}^\top, c \mathbf{A}^\top \right) E[\mathbf{U}^\top \mathbf{U}] \left( \begin{smallmatrix} 0\\ \mathbf{A} \end{smallmatrix} \right) + c E[\varepsilon^\top_X \mathbf{U}] \left( \begin{smallmatrix} 0\\ \mathbf{A} \end{smallmatrix} \right) + \underbrace{E[\varepsilon^\top_Y \mathbf{U}]}_{=0}\left( \begin{smallmatrix} 0\\ \mathbf{A} \end{smallmatrix} \right) \nonumber\\
& + \left( \mathbf{B}^\top, c \mathbf{A}^\top \right) E[\mathbf{U}^\top \varepsilon_X] + c \underbrace{E[\varepsilon_X^\top \varepsilon_X]}_{\sigma_X^2} + \underbrace{E[\varepsilon_Y^\top \varepsilon_X]}_{=0},
\end{align}
where $E[\varepsilon^\top_Y \mathbf{U}]=E[\varepsilon_Y^\top \varepsilon_X]=0$ because $\varepsilon_Y$ is i.i.d. and therefore independent of processes from the past. Note, that the suppressed subscript of $\varepsilon_X$ is $t-\tau$ for $\tau>0$.
Further, $E[Y^\top \mathbf{U}]$ becomes
\begin{align}
E[Y^\top \mathbf{U}] &= \left( \mathbf{B}^\top, c \mathbf{A}^\top \right) E[\mathbf{U}^\top \mathbf{U}] + c E[\varepsilon_X^\top \mathbf{U}] + \underbrace{E[\varepsilon^\top_Y \mathbf{U}]}_{=0},
\end{align}
and
\begin{align}
E[\mathbf{U}^\top X] &=E[\mathbf{U}^\top \mathbf{U}] \left( \begin{smallmatrix} 0\\ \mathbf{A} \end{smallmatrix} \right) + E[\mathbf{U}^\top X].
\end{align}
Then
\begin{align}
&E[Y^\top \mathbf{U}] (E[\mathbf{U}^\top \mathbf{U}])^{-1} E[\mathbf{U}^\top X]=\nonumber\\
&\left( \mathbf{B}^\top, c \mathbf{A}^\top \right) E[\mathbf{U}^\top \mathbf{U}] \left( \begin{smallmatrix} 0\\ \mathbf{A} \end{smallmatrix} \right) + c E[\varepsilon^\top_X \mathbf{U}] \left( \begin{smallmatrix} 0\\ \mathbf{A} \end{smallmatrix} \right) +\nonumber\\
&+ \left( \mathbf{B}^\top, c \mathbf{A}^\top \right) E[\mathbf{U}^\top \varepsilon_X] +\nonumber\\
&+c E[\varepsilon_X^\top \mathbf{U}] (E[\mathbf{U}^\top \mathbf{U}])^{-1} E[\mathbf{U}^\top \varepsilon_X].
\end{align}
Thus, many terms in $E[Y_{\mathbf{U}} ^\top X_{\mathbf{U}} ]$ cancel, and it remains
\begin{align}
&E[Y_{\mathbf{U}} ^\top X_{\mathbf{U}} ] = c\sigma_X^2+ \nonumber\\
&-c\underbrace{E[\varepsilon_X^\top \mathbf{U}] (E[\mathbf{U}^\top \mathbf{U}])^{-1} E[\mathbf{U}^\top \varepsilon_X]}_{(\star)}.
\end{align}
Treating the inverse covariance in the $(\star)$-term with the matrix inversion lemma analogous to Eq.~(\ref{eq:inversion_lemma}) and noting that
\begin{align}
E[\varepsilon_X^\top \mathbf{U}] = \left( 0, \varepsilon_X^\top \mathbf{W} \right),
\end{align}
because $\varepsilon_X$ is independent from the parents $\mathbf{Z}$ of $X$, the $(\star)$-term becomes
\begin{align}
(\star) &= E[\varepsilon_X^\top\mathbf{W}] S_Z^{-1}  E[\mathbf{W}^\top \varepsilon_X].
\end{align}
$S_Z^{-1}$ is again the inverted $N_Y\times N_Y$ matrix of the conditional variance of $\mathbf{W}$ given $\mathbf{Z}$,
\begin{align}
S_Z=  E[\mathbf{W}^\top \mathbf{W}] - E[\mathbf{W}^\top \mathbf{Z}] (E[\mathbf{Z}^\top\mathbf{Z}])^{-1} E[\mathbf{Z}^\top \mathbf{W}].
\end{align}
Along the same derivation the variances are evaluated. All together, the covariances and variances are simplified to
\begin{align} 
E[Y_{\mathbf{U}} ^\top X_{\mathbf{U}} ] &= c \sigma_X^2  \nonumber\\
               &~~~- c E[\varepsilon_X^\top\mathbf{W}] S_Z^{-1}   E[\mathbf{W}^\top \varepsilon_X] \\
E[Y_{\mathbf{U}} ^\top Y_{\mathbf{U}} ] &= \sigma_Y^2 + c^2 \sigma_X^2  \nonumber\\
    &~~~- c^2 E[\varepsilon_X^\top\mathbf{W}] S_Z^{-1}   E[\mathbf{W}^\top \varepsilon_X]\\
E[X_{\mathbf{U}} ^\top X_{\mathbf{U}} ] &= \sigma_X^2  \nonumber\\
               &~~~- E[\varepsilon_X^\top\mathbf{W}] S_Z^{-1}   E[\mathbf{W}^\top \varepsilon_X].
\end{align}

The ``sidepath'' contribution $E[\varepsilon_X^\top\mathbf{W}]$ can be further analysed as follows.
Inserting $t$ and $\tau$ again, the entries of the vector $E[\varepsilon_X^\top \mathbf{W}]$ can be written as
\begin{align}
(E[\varepsilon_X^\top \mathbf{W}])_i = E[\varepsilon_{X,t-\tau} W_{t-g_i}^i],
\end{align}
A simple case where $E[\varepsilon_X^\top \mathbf{W}]$ is zero is given if $\forall i~\tau<g_i$, i.e., all parents of $Y$ are in the past of $X$. But it is interesting to further analyse more complex cases for $\tau \geq g_i$ for any $i$. 
Consider
\begin{align}
E[\varepsilon_{X,t-\tau} W_{t-g_i}^i] &=E[W^i_{t+\tau-g_i} \varepsilon_{X,t} ] \nonumber \\
&= \underbrace{E[W^i_{t+\tau-g_i}X_t]}_{\Gamma_{W_iX}(\tau-g_i)} - \nonumber \\
&~~~ - \sum^{N_X}_{j=1} \Phi_{XZ_j}(h_j) \underbrace{E[W^i_{t+\tau-g_i}Z^j_{t-h_j}]}_{\Gamma_{W_iZ_j}(\tau+h_j-g_i)}.
\end{align}
Analyzing $\Gamma_{W_iX}(\tau-g_i)$,
\begin{align}
&\Gamma_{W_iX}(\tau-g_i) =  \nonumber\\
& =\sum^\infty_{n=0} \sum^N_{r=1} \sum^N_{s=1} \Psi_{W_i r}(n{+}\tau{-}g_i) \Sigma_{rs} \Psi_{Xs}(n),
\end{align}
the linear combination of paths in $\Psi_{Xs}(n)$ can be separated as they either all go through the parents of $X$ or are emanating from $X$, i.e., are of length $n=0$:
\begin{align}
\Psi_{Xs}(n) &= \delta_{X,s} \delta_{n,0} + \sum^{N_X}_{j=1} \Phi_{XZ_j}(h_j) \Psi_{Z_js}(n-h_j)
\end{align}
resulting in 
\begin{align}
&\Gamma_{W_iX}(\tau-g_i) = \nonumber\\
& = \sum^\infty_{n=0} \sum^N_{r=1} \sum^N_{s=1} \Psi_{W_i r}(n{+}\tau{-}g_i) \Sigma_{rs} \delta_{X,s} \delta_{n,0} + \nonumber\\
&~~+ \sum^\infty_{n=0} \sum^N_{r=1} \sum^N_{s=1} \Psi_{W_i r}(n{+}\tau{-}g_i) \Sigma_{rs} \sum^{N_X}_{j=1} \Phi_{XZ_j}(h_j) \Psi_{Z_js}(n{-}h_j)\\
& =  \sum^N_{r=1} \Psi_{W_i r}(\tau{-}g_i) \Sigma_{rX}  + \nonumber\\
&~~+\sum^{N_X}_{j=1} \Phi_{XZ_j}(h_j) \sum^\infty_{n=0} \sum^N_{r=1} \sum^N_{s=1} \underbrace{\Psi_{W_i r}(n{+}\tau{-}g_i) \Sigma_{rs}  \Psi_{Z_js}(n{-}h_j)}_{\Psi_{W_i r}(n{+}\tau{-}g_i{+}h_j) \Sigma_{rs}  \Psi_{Z_js}(n)}\\
&= \sum^N_{r=1} \Psi_{W_i r}(\tau{-}g_i) \Sigma_{rX}   + \sum^{N_X}_{j=1} \Phi_{XZ_j}(h_j) \Gamma_{W_iZ_j}(\tau{+}h_j{-}g_i)
\end{align}
and thus
\begin{align}
(E[\varepsilon_X^\top \mathbf{W}])_i = \sum^N_{r=1} \Psi_{W_ir}(\tau-g_i) \Sigma_{rX}.
\end{align}
$(E[\varepsilon_X^\top \mathbf{W}])_i$ is the linear combination of all paths of length $\tau-g_i$ emanating from $X_t$ or $\mathbf{X}_{r,t}$ with $\Sigma_{rX} \neq 0$ to $W^i_{t+\tau-g_i}$. 

For $\tau<g_i$, $\Psi(n<0)\equiv 0$ and thus for all $i$ $(E[\varepsilon_X^\top \mathbf{W}])_i=0$, confirming the first part of the theorem.
But for all $i$ with $\tau \geq g_i$, $(E[\varepsilon_X^\top \mathbf{W}])_i$ can still be zero if there are no such paths. If that holds for all $i$, the vector $E[\varepsilon_X^\top \mathbf{W}]$ is zero and the simple expression for MIT is obtained.

The MIT for the sidepath example is derived as follows. In this example we have for all $i\neq k:~(E[\varepsilon_X^\top \mathbf{W}])_i=0$ and also assume that additionally for all $i \neq X:~\Sigma_{iX}=0$. Then $\Psi_{W_kX}(1)=d$ and
\begin{align}
(E[\varepsilon_X^\top \mathbf{W}])_k =  \Psi_{W_kX}(1) \Sigma_{XX} = d \sigma_X^2
\end{align}
and with $E[W_k^\top W_k] =d^2 A^\top E[\mathbf{Z}^\top \mathbf{Z}]+d^2 \sigma_X^2 + \sigma^2_{W_k}$ and $E[W_k^\top \mathbf{Z}]=d A^\top E[\mathbf{Z}^\top \mathbf{Z}]$ the conditional variance of $W_k$ is
\begin{align}
S_Z =  d^2 \sigma_X^2 +\sigma^2_{W_k},
\end{align}
and therefore, since $S_Z$ is a scalar,
\begin{align}
E[\varepsilon_X^\top \mathbf{W}] S_Z^{-1}E[\mathbf{W}^\top \varepsilon_Y] = \frac{d^2 \sigma^4_X}{d^2 \sigma_X^2 +\sigma^2_{W_k}},
\end{align}
from which the sidepath MIT follows.


\end{document}